\documentclass[11pt]{amsart}

\newtheorem{cotmb}{Corollary}
\newtheorem{prop}{Proposition}
\newcommand{\cal}{\mathcal}

\newtheorem{tw}{Theorem}
\newtheorem{lemma}{Lemma}

\title{Eta-Einstein condition on twistor spaces of odd-dimensional Riemannian manifolds}
\author{Johann Davidov}
\address{Institute of Mathematics and Informatics \\
Bulgarian Academy of Sciences\\ Acad. G.Bonchev Str. Bl.8\\
1113 Sofia\\ Bulgaria} \email{jtd@math.bas.bg}

\begin{document}

\begin{abstract}
In this note, we find the conditions on an odd-dimensional Riemannian manifolds under which its twistor space is
eta-Einstein.

\vspace{0,1cm} \noindent 2000 {\it Mathematics Subject Classification}. 53C28, 53C25

\vspace{0,1cm} \noindent {\it Key words: twistor spaces, eta-Einstein almost contact metric structures}

\end{abstract}

\thispagestyle{empty}

\maketitle \vspace{0.5cm}

\section{Introduction}

  The construction of the twistor space of an odd-dimensional Riemannian manifold (as well as the one in the
even-dimensional case) can be traced back to R. Penrose \cite{Pen75,Pen83}. These
spaces have been studied by many people, mainly from the point of view of the
$CR$-geometry (see, for example, the literature quoted in \cite{D02,D03,D04}).

  It is convenient to modify slightly the construction in \cite{LB84,R} and to define
the twistor space of an odd-dimensional Riemannian manifold $(M,g)$ as the bundle
${\cal C}$ over $M$ whose fibre at a point $p\in M$ consists of all pairs
$(\varphi,\xi)$ of a skew-symmetric endomorphism $\varphi$ and a unit vector $\xi$ of
$T_pM$ such that $(\varphi,\xi,g_p)$ satisfies the (algebraic) identities in the
definition of an almost contact metric structure (\cite{D02,D03}). The smooth manifold
${\cal C}$ admits two natural partially complex structures ($f$-structures) $\Phi_1$
and $\Phi_2$, a $1$-parameter family of Riemannian metrics $h_t$, $t>0$, compatible
with $\Phi_1$ and $\Phi_2$, and a globally defined $h_t$-unit vector field $\chi$ such
that $(\Phi_{\alpha},\chi,h_t)$, $\alpha=1,2$, is an almost contact metric structure on
${\cal C}$ (we refer to \cite{Blair} for general facts about such structures).
According to \cite{D04} the metrics $h_t$ on ${\cal C}$ are never Einstein. In this
note we consider on ${\cal C}$ a generalization of the Einstein condition adapted for
almost contact metric manifolds, namely the so-called eta-Einstein condition. Let us
note that the almost contact metric structures $(\Phi_{\alpha},\chi,h_t)$,
$\alpha=1,2$, are never Sasakian (even they are never contact \cite{D02,D03}). The
eta-Einstein condition for Sasakian manifolds has been recently discussed in
\cite{BGM}.

   Recall that an almost contact metric structure with contact form $\eta$ and associated metric $h$ on an
odd-dimensional manifold $N$ is said to be eta-Einstein if there exist smooth functions
$a$ and $b$ on $N$ such that
\begin{equation}\label{eq 0.0}
Ricci_{h}(X,Y)=ah(X,Y)+b\eta(X)\eta(Y), ~ X,Y\in TN.
\end{equation}
It is clear that the functions $a$ and $b$ are uniquely  determined - denoting by $s$
and $\xi$ the scalar curvature of $h$ and the $h$-dual vector field of $\eta$, we have
$s=a({\it dim}\>N)+b$, $\>$ $Ricci(\xi,\xi)=a+b$. Note that, in contrast to the
Einstein case, the functions $a$ and $b$ are not constant in the general case.

    The main purpose of this paper is to prove the following

\begin{tw}\label {Pr 0.1}
Let $(M,g)$ be a Riemannian manifold of odd dimension $n\geq 3$. Then its twistor space
${\cal C}$ endowed with the metric $h_t$ and the contact form
$\eta_{t}=h_t(\cdot,\chi)$ is eta-Einstein if and only if $n=3$, the manifold $(M,g)$
is of positive constant curvature $\nu$ and $t\nu=\frac{1}{2}$; in this case we have
$a=\displaystyle{ \frac{3\nu}{2}}$ and $b=-\displaystyle{\frac{\nu}{2}}$ where the
functions $a$ and $b$ on ${\cal C}$ are defined by means of \rm {(\ref{eq 0.0}) }.
\end{tw}
    The proof is based on a coordinate-free formula for the Ricci tensor of the metric $h_t$ in terms of the
curvature of the base manifold $(M,g)$ obtained in \cite{D04}.

\vspace{0.1cm}

   Suppose that the base manifold $M$ is oriented. Then its twistor space ${\cal C}$ is the disjoint
union of the open subsets ${\cal C}_{\pm}$ consisting of the points $(\varphi,\xi)$
that yield $\pm$ the orientation of $T_pM$ via the decomposition $T_pM= {\it
Im}\varphi\oplus {\Bbb R}\xi$ in which the vector space ${\it Im}\varphi$ is oriented
by means of the complex structure $\varphi$ on it. These open sets are diffeomorphic by
the map $(\varphi,\xi)\to (\varphi,-\xi)$ which sends $(\Phi_{\alpha},\chi,h_{t})$ to
$(\Phi_{\alpha},-\chi,h_{t})$. If, in addition, $M$ is of dimension $3$, the map
$(\varphi,\xi)\to\xi$ is a diffeomorphism of ${\cal C}_{\pm}$ onto the unit tangent
bundle $T_{1}M$ of $M$ which sends the characteristic vector field $\chi$ to the
standard characteristic vector field on $T_{1}M$ (but neither of the structures
$\Phi_{\alpha}$ is going to the standard partially complex structure of $T_{1}M$); this
map sends also the metric $h_t$ to the dilation of the Sasaki metric by $2t$ in the
vertical directions. Thus Theorem~\ref{Pr 0.1} implies that the Sasaki metric on the
unit tangent bundle of the $3$-sphere endowed with the $+1$-curvature metric is
eta-Einstein. In fact, by a result of S.Tanno \cite{T}, the Sasaki metric on the unit
tangent bundle of any unit sphere $S^m$ is eta-Einstein and a suitable modification of
this metric gives a homogeneous Einstein metric on $T_{1}S^m$.

\section{Preliminaries}

     Let $V$ be a real $n$-dimensional vector space with a metric $g$.
A partially complex structure on $V$ (or f-structure) of rank $2k$ is an endomorphism
$F$ of $V$ of rank $2k$, $0<2k\leq n$, satisfying $F^3 + F = 0$. We shall say that such
a structure $F$ is compatible with the metric $g$ if the endomorphism $F$ is
skew-symmetric with respect to $g$.

     Given a compatible partially complex structure $F$, we have the orthogonal
decomposition $V={\it Im}F\oplus {\it Ker}F$ and $F$ is a complex structure on the vector space ${\it Im}F$
compatible with the restriction of $g$.

    Denote by $F_{k}(V,g)$ the set of all compatible partially complex structures
of rank $2k$ on $(V,g)$. The group $O(V)$ of orthogonal transformations of $V$ acts
transitively on $F_{k}(V,g)$ by conjugation and $F_{k}(V,g)$ can be identified with the
homogeneous space  $O(n)/(U(k)\times O(n-2k)) $; in particular, ${\it
dim}F_{k}(V,g)=2nk-3k^2-k$. By the results of \cite{AP}, the homogeneous manifold
$F_{k}(V,g)$ admits a unique (up to homotety) invariant  K\"ahler-Einstein structure
structure $(h,{\cal J})$. It can be described in the following way (see, for example,
\cite{D04,D02,D03}): Consider $F_{k}(V,g)$ as a (compact) submanifold of the vector
space $so(V)$ of skew-symmetric endomorphisms of $(V,g)$. Then the tangent space of
$F_{k}(V,g)$ at a point $F$ consists of all endomorphisms $Q\in so(V)$ such that $QF^2
+ FQF + F^{2}Q + Q = 0$. Let $G(S,T)= -\displaystyle{\frac{1}{2}}Trace\,ST$ be the
standard metric on the space $so(V)$. Then the metric $h$ and the complex structure
${\cal J}$ of $T_{F}F_{k}(V,g)$ are given by:
$$
h(P,Q)=2G(P,Q)-G(FPF,Q),~{\cal J}Q=FQ-QF+FQF^2
$$
(the complex structure ${\cal J}$ coincides with the one defined in \cite{R}). It is not hard to compute (see
\cite{D04}) that the scalar curvature of $h$ is equal to $\frac{1}{2}(n-k-1)(2nk-3k^2-k)$.

  Now suppose that $V$ is of odd dimension $n=2k+1$.  A (linear) almost contact metric structure on the
Euclidean space $(V,g)$ is a pair $(\varphi,\xi)$ of an endomorphism $\varphi$ and a unit vector $\xi$ of $V$ such
that $\varphi^{2}x=-x+g(x,\xi)\xi$ and $g(\varphi x,\varphi y)=g(x,y)-g(x,\xi)g(y,\xi)$, $x,y\in V$.

    Denote the set of these structures by $C(V,g)$. It is easy to see
(cf. e.g. \cite{Blair}) that if $(\varphi,\xi)\in C(V,g)$, then $\varphi$ is a compatible partially complex
structure of rank $2k$ and $\varphi(\xi)=0$. Conversely, if $\varphi$ is a compatible partially complex structure
of rank $2k$ and $\varphi(\xi)=0$ for a unit vector $\xi$, then $(\varphi,\xi)\in C(V,g)$.

    The set $C(V,g)$ is a compact submanifold of $so(V)\times V$;
its tangent space at a point $\sigma=(\varphi,\xi)$ consists of all pairs $(Q,\varphi
Q(\xi))$ with $Q\in T_{\varphi}F_{k}(V,g)$. The group $O(V)$ of orthogonal
transformations of $V$ acts transitively on $C(V,g)$ in an obvious way and $C(V,g)$ has
the homogeneous representation $C(V,g)=O(2k+1)/(U(k)\times \lbrace 1 \rbrace)$.  Thus
we have an obvious two-fold covering map
$$
C(V,g)=O(2k+1)/(U(k)\times \lbrace 1 \rbrace)\rightarrow F_{k}(V,g)=
O(2k+1)/(U(k)\times \lbrace 1,-1 \rbrace).
$$
In fact this is the projection map $(\varphi,\xi)\to \varphi$. We lift the
K\"ahler-Einstein structure of $F_{k}(V,g)$ to the manifold $C(V,g)$ by means of this
map  and denote the lifted structure again by $(h,{\cal J})$.

\section{The twistor space $({\cal C},h_t)$ and its Ricci tensor}

  First we recall the definition of the twistor space of partially complex structures \cite{R}.

     Let $(M,g)$ be a Riemannian manifold of dimention $n\geq 3$. Denote
by $\pi:{\cal F}_k\to M$ the bundle over $M$ whose fibre at a point $p\in M$ consists of all compatible partially
complex structures of rank $2k$ on the Euclidean space $(T_{p}M,g_p)$. This is the associated bundle
$$
O(M)\times_{O(n)} F_{k}({\Bbb R}^n)
$$
where $O(M)$ denotes the principal bundle of orthonormal frames on $M$.

       As it is usual in the twistor theory, the manifold ${\cal F}_k$ admits
two partially complex structures $\Phi_1$ and $\Phi_2$ defined as follows \cite{R}:
(recall that a partially complex structure on a manifold is an endomorphism $\Phi$ of
its tangent bundle having constant rank and such that $\Phi^3+\Phi=0$): The Levi-Civita
connection on $M$ gives rise to a splitting ${\cal V}\oplus {\cal H}$ of the tangent
bundle of any bundle associated to $O(M)$ into vertical and horizontal parts. The
vertical space ${\cal V}_f$ of ${\cal F}_k$ at a point $f\in {\cal F}_k$ is the tangent
space at $f$ of the fibre through this point and $\Phi_1|{\cal V}_f$ is defined to be
the complex structure ${\cal J}_f$ of the fibre while $\Phi_2|{\cal V}_f$ is defined as
the conjugate complex structure, i.e. $\Phi_{2}|{\cal V}_f=-{\cal J}_f$. The horizontal
space ${\cal H}_f$ is isomorphic via the differential $\pi_{*f}$ to the tangent space
$T_{p}M$, $p=\pi(f)$, and both $\Phi_1|{\cal H}_f$,  $\Phi_2|{\cal H}_f$ are defined to
be the lift to ${\cal H}_f$ of the endomorphism $f$ of $T_{p}M$.

     The metrics of $F_{k}({\Bbb R}^n)$ and $M$ yield a $1$-parameter family
$h_t, t>0,$ of Riemannian metrics on ${\cal F}_k$ such that $h_t|{\cal V}_f$ is $t$ times the metric $h$ of the
fibre through $f$, $h_t|{\cal H}_f= \pi^{*}g$, and the spaces ${\cal V}_f$ and ${\cal H}_f$ are orthogonal. The
endomorphisms $\Phi_1$ and $\Phi_2$ are skew-symmetric with respect to the metrics $h_t$ and the projection
$\pi:({\cal F}_k,h_t)\to (M,g)$ is a Riemannian submersion with totally geodesic fibres (by the Vilms theorem).

    Now  let $(M,g)$ be a Riemannian manifold of odd dimension $n=2k+1$, $k\geq1$. Slightly modify the twistor
construction in \cite{LB84,R}, we define the twistor space of $(M,g)$ as the bundle ${\cal C}$ over $M$ whose
fibre at a point $p\in M$ consist of all almost contact metric structures on the Euclidean space $(T_pM,g_p)$,
i.e.
$$
{\cal C}=O(M)\times_{O(2k+1)} C({\Bbb R}^{2k+1}).
$$

 Using the Levi-Civita connection on $M$, we can define on ${\cal C}$ a 1-parameter family $h_t, t>0$, and two
partially complex structures $\Phi_1$, $\Phi_2$ of rank $k^2+3k$, skew-symmetric with respect to $h_t$ in the same
way as we did it for the space ${\cal F}_k$. Define a vector field $\chi$ on ${\cal C}$ by setting
$$\chi_{\sigma}=\xi_{\sigma}^{h},~ \sigma=(\varphi,\xi),$$
where $\xi_{\sigma}^{h}$ is the horizontal lift of $\xi$ at the point $\sigma$. Then $(\Phi_{\alpha},\chi,h_t)$,
$\alpha=1,2$, is an almost contact metric structure on ${\cal C}$. The contact distribution of this structure is
obviously ${\cal V}\oplus \lbrace A\in {\cal H}: A\perp\chi \rbrace$.

    If the manifold $M$ is oriented, then
the twistor space ${\cal C}$ is the disjoint union of the open subsets ${\cal C}_{\pm}$ consisting of the points
$(\varphi,\xi)$ that yield $\pm$ the orientation of $T_pM$ via the decomposition $T_pM= {\it Im}\varphi\oplus
{\Bbb R}\xi$ in which the vector space ${\it Im}\varphi$ is oriented by means of the complex structure $\varphi$
on it. The bundles ${\cal C}_{\pm}$ are isomorphic by the map $(\varphi,\xi)\to (\varphi,-\xi)$ which preserves
the horizontal spaces and sends $(\Phi_{\alpha},\chi,h_{t})$ to $(\Phi_{\alpha},-\chi,h_{t})$.

       The natural covering map $C({\Bbb R}^{2k+1})\to F_{k}({\Bbb R}^{2k+1})$
is $O(2k+1)$-equivariant, so it determines a bundle map
$$
{\cal C}\to {\cal F}_{k},\, (\varphi,\xi)\to \xi,
$$
which is a two-fold covering. This map preserves the vertical and horizontal spaces, the metrics and the partially
complex structures of ${\cal C}$ and ${\cal F}_{k}$. If $M$ is oriented, each of the spaces ${\cal C}_{+}$ and
${\cal C}_{-}$ is isomorphic to ${\cal F}_{k}$.

    The curvature of the Riemannian manifold $({\cal F}_k,h_t)$ has been computed in \cite{D04} by means of the
O'Neill formulas. The computation there immediately gives the curvature of the manifold $({\cal C},h_t)$ since the
map $(\varphi,\xi)\to \xi$ above is a Riemannian covering. To formulate the corresponding result for the curvature
of $({\cal C},h_t)$, we shall introduce some notations. Denote by $A(TM)$ the bundle of skew-symmetric
endomorphisms of $TM$ and consider ${\cal C}$ as a submanifold of the bundle $\pi:A(TM)\oplus TM\to M$. Then the
inclusion of ${\cal C}$ is fibre-preserving and the horizontal subspace of $T_{\sigma}{\cal C}$ at a point
$\sigma=(\varphi,\xi)$ coincides with the horizontal space of $A(TM)\oplus TM$ at that point. The vertical space
of ${\cal C}$ at $\sigma$, considered as a subspace of the vertical space of $A(TM)\oplus TM$ at $\sigma$,
consists of all pairs $(Q,\varphi Q(\xi))$ for which $Q\in A(T_pM)$, $p=\pi(\sigma)$,  and satisfies the identity
$Q\varphi ^{2} + \varphi ^{2}Q+\varphi Q\varphi+Q=0$.  Further we shall freely make use of the standard isometric
identification $A(TM)\cong\Lambda^{2}TM$ that assigns to each $a\in A(T_{p}M)$ the 2-vector $a^{\wedge}$ for which
$g(aX,Y)=g(a^{\wedge},X\wedge Y), X,Y\in T_{p}M$ (the metric on $\Lambda^{2}TM$ is given by $g(X_1\wedge
X_2,X_3\wedge X_4)= g(X_1,X_3)g(X_2,X_4)-g(X_1,X_4)g(X_2,X_3)$). If $Q\in T_{\varphi}F_k(T_pM,g_p)\subset
A(T_pM)$, the element $Q^{\wedge}$ of $\Lambda^{2}T_{p}M$ will be also denoted by $Q$. For brevity, denote by
$m_{\varphi}$ the image of the tangent space $T_{\varphi}F_k(T_pM,g_p)$, $p=\pi(\sigma)$, under the identification
$A(T_pM)\cong\Lambda^2T_pM$. Finally, let ${\cal R}:\Lambda^{2}TM\to\Lambda^{2}TM$ be the curvature operator of
$(M,g)$; it is defined by $g({\cal R}(X_1\wedge X_2),X_3\wedge X_4)= g((\nabla_{[X_1,X_2]}-
[\nabla_{X_1},\nabla_{X_2}])X_3,X_4)$). Now \cite[Proposition 2]{D04} implies the following

\vspace{0.1cm}
\begin{prop}\label {Pr 3.1}
Let $(M,g)$ be a Riemannian manifold of odd dimension $n=2k+1$, $k\geq1$. Then the the
Ricci tensor $c_t$ of its twistor space $({\cal C},h_t)$ is given as follows: For any
$E\in T_{\sigma}{\cal C}$, $\sigma=(\varphi,\xi)$, setting $X=\pi_{\ast}E$, $(Q,\varphi
Q(\xi))={\cal V}E$ (the vertical component of $E$), we have:
$$
\begin{array}{lll}
c_t(E,E)&=&c_M(X,X)-2tTrace(Z\to (\nabla_{Z}R)({\cal
J}Q,X))+\\
                 & &2t^2||{\cal R}({\cal J}Q)||^2_{g}
-2t||\imath_X\circ{\cal R}|m_{\varphi}||^2_{h,g}+\displaystyle{\frac{1}{2}}k||Q||^2_{h}
\end{array}
$$
where $c_M$ is the Ricci tensor of $(M,g)$, $\imath_X:\Lambda^2TM\to TM$ is the interior product and
$||\cdot||_{h,g}$ is the norm of the metric on the space of linear maps $m_{\varphi}\to T_{\pi(\sigma)}M$ induced
by the metrics $h$ on $m_{\varphi}$ and $g$ on $T_{\pi(\sigma)}M$.
\end{prop}

\vspace{0.1cm}
\begin{cotmb}\label {Cor 3.3}
If $(M,g)$ is of constant curvature $\nu$, then the Ricci tensor $c_t$ of $({\cal C},h_t)$ is given by
$$
c_t(E,E)=2k\nu(1-t\nu)||X||^2_{g}+t\nu^2||\varphi X||^2_{g}+
$$
$$
\frac{1}{2}(k+2t^2\nu^2)||Q||^2_{h}+t^2\nu^2h(\varphi\circ Q\circ \varphi,Q)
$$
where $\sigma=(\varphi,\xi)\in {\cal C}$, $E\in T_{\sigma}{\cal C}$, $X=\pi_{\ast}E$ and $(Q,\varphi Q(\xi))={\cal
V}E$.
\end{cotmb}

\section{Proof of the Theorem}

   Suppose that $({\cal C},h_t,\chi)$ is eta-Einstein. Then, by Proposition 1, there exist smooth functions $a$
and $b$ on ${\cal C}$ such that for every point $\sigma=(\varphi,\xi)\in {\cal C}$, every $h$-orthonormal basis
$\{U_{\alpha}\}$ of $T_{\varphi}F_k(T_p,g_p)$, $p=\pi(\sigma)$, and every $X\in T_pM$, $Q\in
T_{\varphi}F_k(T_p,g_p)$ the following two equations are satisfied
\begin{equation}\label{eq 4.1}
c_M(X,X) -2t\sum_{\alpha=1}^{k^2+k}||R(U_{\alpha})X||^2=a(\sigma)||X||^2 + b(\sigma)g(X,\xi)^2
\end{equation}
\begin{equation}\label{eq 4.2}
2t^2||{\cal R}(Q)||^2_g +\frac{1}{2}k||Q||^2_h=a(\sigma)t||Q||^2_h
\end{equation}

\vspace{0.1cm}

\begin{lemma}\label {Descend}
The functions $a$ and $b$ on ${\cal C}$ descend to smooth functions $\bar a$ and $\bar b$ on $M$.
\end{lemma}

\begin{proof}
We have
\begin{equation}\label{eq 4.3}
c_t(E',E'')=ah_t(E',E'')+bh_t(E',\chi)h_t(E'',\chi)
\end{equation}
for every $E',E''\in T{\cal C}$.

 Take a point  $\sigma=(\varphi,\xi)\in {\cal C}$  and set $p=\pi(\sigma)$. Let $e_1,...,e_{2k+1}$ be an
orthonormal basis of $T_pM$ with $e_{2k+1}=\xi$ and let $V_1,...,V_{k^2+k}$ be a $h_t$-orthonormal basis of the
vertical space ${\cal V}_{\sigma}$ of ${\cal C}$. Denote by $\{A_i\}$ the $h_t$-orthonormal basis
$\{e_1^h,...,e_{2k+1}^h,V_1,...,V_{k^2+k}\}$ of $T_{\sigma}{\cal C}$. Denote by $D$ the Levi-Civita connection of
the metric $h_t$ on ${\cal C}$. Then a standard application of the differential Bianchi identity gives:
$$
\sum_{i}A_m(c_t(A_i,A_i))=\sum_{i,j}h_t((D_{A_m}R_t)(A_i,A_j,A_i),A_j)=2\sum_{i}A_i(c_t(A_i,A_m))
$$
where $R_t$ is the curvature tensor of the metric $h_t$ on ${\cal C}$. For every
$l=1,...,2k+1$ and $\beta=1,...,k^2+k$, we have $c_t(e_l^h,V_{\beta})=0$ by (\ref{eq
4.3}). Therefore
$$
(k^2+3k+1)V_{\beta}(a)+V_{\beta}(b)=\sum_{i}V_{\beta}(c_t(A_i,A_i))=2\sum_{\gamma=1}^{k^2+k}V_{\gamma}(c_t(V_{\gamma},V_{\beta}))=
2V_{\beta}(a).
$$
It follows that the function $(k^2+3k-1)a+b$ is constant on the fibers of ${\cal C}$. Denote by $s$ the scalar
curvature of $(M,g)$. Then (\ref{eq 4.1}) and (\ref{eq 4.2}) imply that
$$
s(p)+t^{-1}(k^3+k^2)=(2k^2+4k+1)a(\sigma)+b(\sigma).
$$
Thus, the function $(2k^2+4k+1)a+b$ is also constant on the fibers of ${\cal C}$. This proves the lemma.
\end{proof}

\vspace{0.1cm}

  Let $p\in M$ and let $e_1,e_2,....,e_n$ be an orthonormal basis of $T_pM$, $n=2k+1$.
Set $e_{ij}=e_i\wedge e_j$. Assume that $n=2k+1\geq 5$.

\vspace{0.1cm}
\begin{lemma}\label {L}
For any $X\in T_pM$ we have
\begin{equation}\label{eq 4.5}
\begin{array}{c}
c_M(X,X) -t\sum_{p=1}^{k}\sum_{j=2p+1}^{n}(||R(e_{2p-1,j})X||^2 + ||R(e_{2p,j})X||^2)= \\
                                                                                        \\
\bar a(p)||X||^2 +\bar b(p)g(X,e_n)^2.
\end{array}
\end{equation}
\end{lemma}

\begin{proof}
Consider the point $\sigma=(\varphi,\xi)\in{\cal C}$ defined by $\varphi=e_{12}+....+e_{2k-1,2k}$ and
$\xi=e_{2k+1}$. Set
$$
\begin{array}{c}
A_{pq}^{'}=\frac{1}{\sqrt 2}(e_{2p-1,2q-1}-e_{2p,2q}), ~ A_{pq}^{''}=\frac{1}{\sqrt 2}(e_{2p-1,2q}+e_{2p,2q-1})\\
                                                                                                             \\
B_{r}^{'}=\frac{1}{\sqrt 2}e_{2r-1,2k+1},  ~ B_{r}^{''}=\frac{1}{\sqrt 2}e_{2r,2k+1}
\end{array}
$$
where $p=1,...,k-1$, $q=p+1,...,n$, $r=1,...,k$. Then $\{U_{\alpha}\}=
\{A_{pq}^{'},A_{pq}^{''},B_{rs}^{'},B_{rs}^{''}\}$ is a $h$-orthonormal basis of
$T_{\varphi}F_k(T_pM,g_p)$ (such that ${\cal J}A_{pq}^{'}=A_{pq}^{''}$, ${\cal
J}B_{rs}^{'}=B_{rs}^{''}$). Writing (\ref{eq 4.1}) for this basis we get an identity
which involves the basis $e=(e_1,...,e_n)$; in view of Lemma 1, its right-hand side
depends on the choice of the vector $e_n$ and does not depend on the particular choice
of $\varphi$. We denote by $(\ref{eq 4.1})_e$ the identity we obtain in this way. Set
$e^{'}=(-e_1,e_2, ...,e_n)$. Then the identity $(\ref{eq 4.1})_e - (\ref{eq
4.1})_{e^{'}}$ reads as
\begin{equation}\label{eq 4.0}
\begin{array}{c}
g(R(e_{13})X,R(e_{24})X)-g(R(e_{14})X,R(e_{23})X)+......+  \\
g(R(e_{1,2k-1})X,R(e_{2,2k})X)-g(R(e_{1,2k})X,R(e_{2,2k-1})X)=0.
\end{array}
\end{equation}
Now we  apply (\ref{eq 4.0}) for the bases $e=(e_1,...,e_4,...,e_n)$ and $e^{''}=(e_1,...,-e_4,...,e_n)$. Then the
identity $(\ref{eq 4.0})_e - (\ref{eq 4.0})_{e^{''}}$ gives
\begin{equation}\label{eq 4.4}
g(R(e_{13})X,R(e_{24})X)-g(R(e_{14})X,R(e_{23})X)=0.
\end{equation}
It follows that $g(R(e_{2p-1,2q-1})X,R(e_{2p,2q})X)=g(R(e_{2p-1,2q})X,R(e_{2p,2q-1})X)$
for $p=1,2,...,k-1$, $q=p+1,p+2,...,n$. The latter identity, $(\ref{eq 4.1})_e$ and
Lemma 1 imply Lemma 2.
\end{proof}

\vspace{0.1cm}

  The proofs of the next two lemmas go in the same lines as the proofs of Lemmas 6 and 7 in \cite{D04}
(in view of Lemmas 1 and 2 above) and will be omitted.

\vspace{0.1cm}
\begin{lemma}\label {L4}
If $i,j,,l,m\in \{1,...,n\}$ are four different indexes, then
$$
g(R(e_{ij})X,R(e_{lm})X)=0
$$
for any $X\in T_pM$.
\end{lemma}

\vspace{0.1cm}

\begin{lemma}\label {LN}
For any $i\neq j$, $l\neq m$ and $X\in T_pM$, we have
$$
||R(e_{ij})X||=||R(e_{lm})X||.
$$
\end{lemma}

\vspace{0.1cm}

\newpage

\noindent {\it Proof of the theorem in the case $n=2k+1\geq 5$}.

  Assume that $({\cal C},h_t,\chi)$ is eta-Einstein. Let $p\in M$ and let $e_1,....,e_n$, $n=2k+1$, be an orthonormal basis
of $T_pM$. Set $\varphi=e_{12}+....+e_{2k-1,2k}$ and $\xi=e_{2k+1}$. Set also
$$
Q_1=e_{13}-e_{24}, \> Q_2=e_{1,2k+1}.
$$
Then $Q_1,Q_2\in T_{\varphi}F_k(T_pM,g_p)$ and according to (\ref{eq 4.2}) and Lemma 1
$$
||{\cal R}(Q_1)||^2=||{\cal R}(Q_2)||^2=t^{-2}[\bar a(p)t-\frac{1}{2}k].
$$
On the other hand, by Lemmas 3 and 4
$$
||R(Q_1)X||^2=2||R(e_{13})X||^2=2||R(e_{1,2k+1})X||^2=2||R(Q_2)X||^2
$$
for any $X\in T_pM$. Therefore $||{\cal R}(Q_1)||^2=2||{\cal R}(Q_2)||^2$ and we see that ${\cal R}(Q_1)={\cal
R}(Q_2)=0$. Then
$$
||{\cal R}(e_{13})||=2^{-1/2}||{\cal R}(Q_1)||=0 ~ \mbox { and } ~ \bar a(p)t-\frac{1}{2}k=0.
$$
The first of these identities implies $R=0$. Now taking $X=e_1$ in (\ref{eq 4.1}) we see that $\bar a(p)=0$ which
contradicts to the second identity above.

\vspace{0.1cm}

\noindent {\it Proof of the theorem in the case $n=3$}.

    Let $c_M$ be the Ricci tensor of $(M,g)$ and $\rho:T_pM\to T_pM$ the symmetric operator corresponding to $c_M$..
Denote by $s$ the scalar curvature of $(M,g)$. It is well-known that the curvature operator of a $3$-dimensional
Riemannian manifold is given by
\begin{equation}\label{eq CO}
{\cal R}(X\wedge Y)=-\frac{s}{2}X\wedge Y +\rho(X)\wedge Y + X\wedge\rho(Y)
\end{equation}
for $X,Y\in TM$ (see e.g. \cite[Sec. 1 G]{Be}).
   Let $p\in M$ and put
$$
\lambda=\frac{1}{2t^2}[\bar a(p)t-\frac{1}{2}k].
$$
Let $e_1,e_2,e_3$ be an arbitrary orthonormal basis of $T_pM$. Consider the point $\sigma=(\varphi,\xi)\in {\cal
C}$ with $\varphi=e_{12}$, $\xi=e_3$. Then $e_{13}$, $e_{23}\in T_{\varphi}F_1(T_pM,g_p)$ and by (\ref{eq 4.2}) we
have
$$
||{\cal R}(e_{13})||^2=||{\cal R}(e_{23})||^2=2\lambda ~ \mbox { and } ~ g({\cal R}(e_{13}),{\cal R}(e_{23}))=0.
$$
It follows that either $\lambda=0$  or the operator $T=(1/\sqrt{2\lambda}){\cal R}$ is orthogonal. If $\lambda=0$,
then ${\cal R}=0$  and identity (\ref{eq 4.1}) implies $\bar a(p)=0$. This together with $\lambda=0$ gives $k=0$,
a contradiction. Thus, the operator $T$ is orthogonal. Since this operator is also symmetric, its square is equal
to $Id$. Therefore the eigenvalues of $T$ are $+1$ or $-1$. Suppose that both $+1$ and $-1$ are eigenvalues of the
operator $T$ and denote by $\alpha$ and $\beta$ the dimensions of the corresponding eigenspaces. Then
\begin{equation}\label{eq Scal}
\frac{s}{2}=Trace\>{\cal R}=(\alpha-\beta)\sqrt{2\lambda}.
\end{equation}
Further, by (\ref{eq 4.1}) and (\ref{eq 4.2}), we have that
$$
s-8t\lambda=3\bar a+\bar b,
$$
therefore
\begin{equation}\label{eq Lambda}
\lambda=\frac{2ts-3-2\bar{b}t}{28t^2}.
\end{equation}

 Set $c_{ij}=c_M(e_i,e_j)$.  Now take the point $\sigma=(e_{13},e_2)\in {\cal C}$ and apply (\ref{eq 4.1}) with $U_1=
\frac{1}{\sqrt 2}e_{12}$, $U_2=\frac{1}{\sqrt 2}e_{23}$, $X=e_1$. This gives
\begin{equation}\label{eq C1}
c_{11}-t(||R(e_{12})e_1||^2 +||R(e_{23})e_1||^2)=\bar a(p).
\end{equation}
Similarly, considering the point $(e_{23},e_1)\in {\cal C}$ and applying (\ref{eq 4.1}), we get
$$
c_{11}-t(||R(e_{12})e_1||^2 + ||R(e_{13})e_1||^2)=\bar a(p)+\bar b(p).
$$
Hence, $\bar b(p)=t(||R(e_{23})e_1||^2 - ||R(e_{13})e_1||^2)$, which, in view of (\ref{eq CO}), implies
$$
\bar b(p)=[c_{12}^2-(\frac{s}{2}-c_{22})^2].
$$
Similarly, we have also
$$
\bar b(p)=[c_{23}^2-(\frac{s}{2}-c_{33})^2] ~ \mbox { and } ~ \bar b(p)=[c_{13}^2-(\frac{s}{2}-c_{11})^2].
$$
Adding the last three equalities and setting $\mu=c_{12}^2+c_{23}^2+c_{13}^2$, we obtain
\begin{equation}\label{eq b}
3\bar b(p)=t[\mu-(\frac{s}{2}-c_{11})^2-(\frac{s}{2}-c_{22})^2-(\frac{s}{2}-c_{33})^2].
\end{equation}
Next, an obvious application of (\ref{eq 4.1}) for the point $(e_{12},e_3)\in {\cal C}$, gives
$$
c_{11}-t(||R(e_{13})e_1||^2 + ||R(e_{23})e_1||^2)=\bar a(p).
$$
From the latter identity and (\ref{eq C1}) we see that
$$
||R(e_{12})e_1||=||R(e_{13})e_1||.
$$
This and (\ref{eq CO}) imply
$$
(-\frac{s}{2}+c_{11}+c_{22})^2=(-\frac{s}{2}+c_{11}+c_{33})^2
$$
which gives $c_{11}(c_{22}-c_{33})=0$. Similarly, $c_{22}(c_{33}-c_{11})=0$ and
$c_{33}(c_{11}-c_{22})=0$. It follows that either $c_{11}=c_{22}=c_{33}$ or two of the
numbers $c_{11}$, $c_{22}$, $c_{33}$ are equal to 0 and the third one is different from
zero.

$ 1)$. Suppose $c_{11}=c_{22}=c_{33}$. Then, by (\ref{eq b}),
$$
3\bar b(p)= t\mu -\frac{ts^2}{12}.
$$
Now it follows from (\ref{eq Scal}) and (\ref{eq Lambda}) that $ts$ satisfies the equation
\begin{equation}\label{eq ts}
(7-\frac{m}{9})x^2-4mx +6m + \frac{4mt^2\mu}{3}=0
\end{equation}
where $m=(\alpha-\beta)^2$. This fact implies
$$
4m^2\geq (7-\frac{m}{9})(6m+\frac{4mt^2\mu}{3})\geq (7-\frac{m}{9})6m.
$$
Thus we see that $m\geq 9$. On the other hand $|\alpha-\beta|<dim\>\Lambda^2 T_pM=3$, so $m<9$, a contradiction.

 $2)$. Assume that $c_{11}=c_{22}=0$. In this case, according to (\ref{eq b}), we have
$$
3\bar b(p)= t\mu -\frac{3ts^2}{4}
$$
and, in view of (\ref{eq Scal}) and (\ref{eq Lambda}), $ts$ satisfies the equation
$$
6x^2 - 4mx +6m + \frac{4mt^2\mu}{3}=0.
$$
This implies $m\geq 9$ and we come again to a contradiction.

  It follows that  either $T=Id$ or $T=-Id$ on the whole space
$\Lambda^2T_pM$. Therefore the sectional curvature of $M$ at each point $p$ is constant and the classical Schur
theorem implies that $M$ is of constant curvature, say $\nu$.  Moreover $ts$ satisfies equation (\ref{eq ts}) with
$m=9$ and $\mu=0$. Thus  $ts=3$, i.e. $t\nu=\frac{1}{2}$.

   Conversely, suppose that $M$ is a $3$-dimensional Riemannian manifold of positive constant curvature $\nu$ and
take $t=\frac{1}{2\nu}$. Let $\sigma=(\varphi,\xi)\in {\cal C}$, $E\in T_{\sigma}{\cal C}$, $X=\pi_{\ast}E$ and
${\cal V}E=(Q,\varphi Q(\xi))$. Since $dim\>M=3$, we have $\varphi\circ Q\circ\varphi=0$ and Corollary 1 gives
$$
\begin{array}{c}
c_t(E,E)=(2\nu-t\nu^2)||X||^2 +\displaystyle{\frac{1}{2t}}(1+2t^2\nu^2)||Q||_{h_t}^2 - t\nu^2g(X,\xi)^2=\\
                                                                                          \\
         \displaystyle{ \frac{3\nu}{2}}||E||_{h_t}^2 -\displaystyle{\frac{\nu}{2}}g(X,\xi)^2.
\end{array}
$$

\end{document}